\documentclass[reqno]{amsart}
\usepackage{amssymb}

\newtheorem{theorem}{Theorem}[section]
\newtheorem{proposition}[theorem]{Proposition}
\newtheorem{lemma}[theorem]{Lemma}

\newtheorem{definition}[theorem]{Definition}

\theoremstyle{remark}
\newtheorem{remark}[theorem]{Remark}

\numberwithin{equation}{section}

\begin{document}

\title[Pieri Formulas]
{Pieri Formulas for Macdonald's Spherical Functions and Polynomials}

\author{J.F.  van Diejen}
\author{E. Emsiz}
\address{
Instituto de Matem\'atica y F\'{\i}sica, Universidad de Talca,
Casilla 747, Talca, Chile}

\subjclass[2000]{Primary: 05E05; Secondary 33D52, 17B22, 20F55} \keywords{symmetric functions, Hall-Littlewood polynomials, Macdonald polynomials, Pieri formulas, root systems, reflection groups}

\thanks{Work supported in part by the {\em Fondo Nacional de Desarrollo
Cient\'{\i}fico y Tecnol\'ogico (FONDECYT)} Grants \# 1090118 and
3080006, and by the {\em Programa Reticulados y
Ecuaciones} of the Universidad de Talca.}

\date{April 2010}

\begin{abstract}
We present explicit Pieri formulas for Macdonald's spherical functions (or generalized Hall-Littlewood polynomials associated with root systems) and their $q$-deformation the Macdonald polynomials. For the root systems of type $A$, our Pieri formulas recover the well-known Pieri formulas for the Hall-Littlewood and Macdonald symmetric functions
due to Morris and Macdonald as special cases.
\end{abstract}

\maketitle

\section{Introduction}\label{sec1}
Given a system of generators and a basis for a commutative ring, the corresponding Pieri formulas describe the action of multiplication by the generators in terms of the basis and thus
completely encode the multiplicative structure of the ring in question.
Formulas of this type constitute a fundamental tool in algebraic combinatorics, where they are used to study and proof properties of special bases of rings of symmetric functions, such as e.g. the Schur functions, the Hall-Littlewood polynomials, and the Macdonald polynomials \cite{mac:symmetric}.
Specifically, the Pieri formulas for the Macdonald polynomials give rise to a straightforward proof of Macdonald's normalization and principal specialization formulas \cite{mac:symmetric} and have, more recently, also been instrumental for obtaining closed expressions for the expansion coefficients of these polynomials in the basis of elementary symmetric functions \cite{las-sch:inversion}. More classically, the Pieri formulas for the Schur functions have proven to be crucial in Schubert calculus for establishing the link with the cohomology ring of Grassmannian varieties permitting, for instance, the computation of intersection numbers of Schubert classes 
\cite{man:symmetric} (cf. also  \cite{lam-lap-mor-sch:affine} for recent developments concerning analogous links between the Schubert calculus on affine Grassmannians and $k$-Schur functions with a prominent role played again by the relevant Pieri formulas).

The purpose of this paper is to present an explicit Pieri formula for Macdonald's zonal spherical functions (on $p$-adic symmetric spaces) known also by the name generalized Hall-Littlewood polynomials associated with root systems \cite{mac:spherical,mac:orthogonal,nel-ram:kostka}.
For the root system of type $A_n$, corresponding to the group $SL(n+1)$, our Pieri formula
contains the Pieri formulas for the Hall-Littlewood polynomials found by Morris  \cite{mor:note} as special cases. More generally, for the classical root systems our Pieri formula specializes to a complete system of Pieri formulas (i.e. corresponding to a complete system of generators) whereas for the exceptional root systems only part of a generating system is covered. The main idea of the proof of our Pieri formula is by degeneration from a more general Pieri formula for the $q$-deformation of Macdonald's spherical functions: the Macdonald polynomials  \cite{mac:orthogonal,mac:affine}.

Recently the study of Hall-Littlewood polynomials and, more generally,  Macdonald spherical functions has experienced some remarkable boosts prompting our demand for explicit Pieri formulas. Important developments worth emphasizing in this context are the combinatorial formulas
for these polynomials \cite{sch:galeries,ram:alcove,len:hall-littlewood,len:haglund}
and for the corresponding Littlewood-Richardson coefficients or structure constants
\cite{par:buildings,sch:galeries,ram:alcove}, and also their linking with the Schubert calculus of isotropic Grassmannians \cite{tam:giambelli}. As an application along very different lines, the Pieri formulas in the present paper also give rise to explicit formulas for the commuting quantum integrals of certain discrete integrable many-particle systems whose scattering behaviour was studied in Ref. \cite{die:scattering}.

The paper is organized as follows. In Section \ref{sec2} we present our Pieri formula.
The proof is relegated to Sections \ref{sec3}, \ref{sec4} and \ref{sec5}. Specifically, first it is shown how the Pieri formula for the Macdonald spherical functions arises as a ($q\to 0$) limiting case of an analogous Pieri formula for the Macdonald polynomials \cite{mac:orthogonal,mac:affine} (Section \ref{sec3}). The latter Pieri formula follows in turn from a difference equation for the Macdonald polynomials upon invoking the duality symmetry (Section \ref{sec4}). Finally, the difference equation for the Macdonald polynomials is proved using a residue calculus for root systems (Section \ref{sec5}). Some useful properties of a special class of small dominant weights appearing frequently throughout the paper have been collected in a short appendix.


\section{Pieri formula for Macdonald spherical functions}\label{sec2}

Let $R$ be an irreducible reduced crystallographic root system spanning a real (finite-dimensional) Euclidean vector 
space $E$ with inner product $\langle \cdot ,\cdot \rangle$.
Following standard conventions \cite{bou:groupes}, let us write $Q$, $P$, and $W$, respectively, for the root lattice, the weight lattice, and the Weyl group associated with $R$. The semigroup of the root lattice generated by a (fixed) choice of positive roots $R^+$
is denoted by $Q^+$ and we write $P^+$  for the corresponding cone of dominant weights.
The standard basis for the group algebra
$\mathbb{C}[P]$ is given by the formal exponentials
$e^\lambda$, $\lambda\in P$ characterized by the relations $e^0=1$, 
$e^\lambda e^\mu=e^{\lambda +\mu}$. For any fixed $q\in (0,1)$, the elements of the group algebra can be thought of as functions on $E$ through the evaluation homomorphism
$e^\lambda (x):= q^{\langle \lambda ,x\rangle}$,  $x\in E$.
The Weyl group acts on $\mathbb{C}[P]$ 
via $we^\lambda :=e^{w\lambda}$, $w\in W$.

Let
$\rho_{t;R}:=\frac{1}{2} \sum_{\alpha\in R^+}\log_q (t_\alpha)\alpha$. Here
$t:R\cup R^\vee \rightarrow (0,1)$ (with $R^\vee:=\{ \alpha^\vee \mid \alpha\in R\}$
where $\alpha^\vee:=2\alpha/\langle \alpha,\alpha\rangle$)
denotes a root multiplicity function
such that  $t_{w\alpha}=t_\alpha$ for all $w\in W$ and $t_{\alpha^\vee}=t_\alpha$ for all $\alpha\in R$. The Macdonald spherical functions constitute a special basis of the $W$-invariant subalgebra $\mathbb{C}[P]^W\subset \mathbb{C}[P]$ parametrized by $t$.

\begin{definition}[Macdonald Spherical Function \cite{mac:spherical,mac:orthogonal}]\label{msf:def}
For $\lambda\in P^+$ the Macdonald spherical function is defined by the formula
\begin{subequations}
\begin{equation}\label{hl1}
P_{\lambda} := \frac{e^\lambda(\rho_{t;R^\vee})}{W(t)} \sum_{w\in W}
\Bigl(e^{w\lambda } \prod_{\alpha\in R^+}\frac{1-t_\alpha
e^{-w\alpha }}{1-
e^{-w \alpha}} \Bigr) 
\end{equation}
with
\begin{equation}\label{hl2}
W(t):=\sum_{w\in W} \prod_{\substack{\alpha \in R^+\\  w\alpha\not\in R^+}} t_\alpha
=\prod_{\alpha\in R^+} \frac{1-t_\alpha e^{\alpha}(\rho_{t,R^\vee}) }{1-e^{\alpha}(\rho_{t,R^\vee})  } 
\end{equation}
\end{subequations}
(where the normalization is chosen such that $P_\lambda (\rho_{t;R^\vee})=1$).
\end{definition}

The cone of dominant weights is partially ordered by the dominance ordering: for $\lambda ,\mu\in P^+$,
$\lambda\geq \mu$ if and only if $\lambda-\mu\in Q^+$.
For $\lambda\in P^+$ we distinghuish the highest-weight system
$P_\lambda^+:=\{ \mu\in P^+\mid \mu \leq \lambda \}$, its cardinality
$n_\lambda:= | P_\lambda^+|$, and the monomial symmetric function
$ m_{\lambda ;R} :=\sum_{\nu\in W \lambda} e^\nu $.

\begin{definition}[Small Weights]\label{small:def}
A dominant weight $\omega $ will be called {\em small} if $\langle \omega,\alpha^\vee\rangle \leq 2$ for any positive root  $\alpha$.
\end{definition}

Some relevant properties enjoyed by
the small weights have been collected in Appendix \ref{sw:app} below.
For instance, a particularly important feature to bear in mind is that for $\omega$ small  the
highest-weight system  $P^+_\omega$ forms a linear chain of small weights
(cf.  Lemma \ref{sw-chain:lem}).

For any $x\in E$, we define the stabilizer subgroup
$W_{x}:=\{ w\in W \mid wx=x\} $, its root subsystem
$R_{x}:= \{ \alpha \in R\mid \langle x,\alpha\rangle =0\} $, the corresponding choice of positive roots $R_{x}^+:=R^+\cap R_{x}$, and the shortest element $w_x\in W$ mapping 
$x$ into the (closed) dominant chamber.
With these (standard) notations we are now in the position to formulate our main result.

\begin{theorem}[Pieri Formula for Macdonald Spherical Functions]\label{PF:thm}
Let $\lambda, \omega \in P^+$ with $\omega $ small. Then
\begin{subequations}
\begin{equation}
E_{\omega;R} P_\lambda =
\sum_{\substack{\mu_1<\cdots <\mu_\ell=\omega\\ \ell =1,\ldots ,n_\omega}}
(-1)^{\ell -1} \!\!\!\!\!\!\!\!
\sum_{\substack{\nu_k\in W_{\nu_{k-1}}  (  w^{-1}_{\nu_{k-1}} \mu_k  )   \\ k=1,\ldots ,\ell}} \!\!\!\!\!\!\!\!
  P_{\lambda +\nu_1}
\prod_{1\leq k\leq \ell}
V_{\nu_k;\nu_{k-1}}(\lambda)  ,
\end{equation}
where 
\begin{equation}\label{Eomega}
E_{\omega;R}  :=
\sum_{\substack{\mu_1<\cdots <\mu_\ell=\omega\\ \ell =1,\ldots ,n_\omega}}
(-1)^{\ell -1} m_{\mu_1;R}
\prod_{1\leq k\leq \ell-1}
m_{\mu_{k+1};R_{\mu_{k}}}(\rho_{t;{R^\vee_{\mu_{k}}}}) 
\end{equation}
and
\begin{eqnarray}
&& V_{\nu;\eta}(\lambda) :=
 e^{-\nu}(\rho_{t;R^\vee_{ \eta}}-\rho_{t;R^\vee_{\eta}\cap R^\vee_{\lambda}})   \\
&& \times \prod_{\substack{\alpha\in R_{\eta} \\    \langle \alpha^\vee ,\lambda\rangle \in \{ 0,-1\}  \\  \langle \alpha^\vee ,\lambda+\nu\rangle>0 }}
t_\alpha^{-\langle \nu ,\alpha^\vee\rangle /2}
\left(\frac{1-t_\alpha e^{\alpha}(\rho_{t;R^\vee}) }{1-e^{\alpha}(\rho_{t;R^\vee})  } \right) .
\nonumber
\end{eqnarray}
\end{subequations}
\end{theorem}

In these formulas we have adopted the convention that $R_{\nu_0}:=R$,  $W_{\nu_0}:=W$ and $w_{\nu_0}:=\text{Id}$, so
$V_{\nu_1;\nu_{0}}=V_{\nu_1; 0}=V_{\nu_1; -}$ (where the dash indicates that the second argument is absent). More generally, throughout the text
indexed weights like $\nu_k$ and $\mu_k$ are understood to be {\em absent} when the index $k$ is not positive.

It seems a daunting task to verify the above Pieri formula directly. In the next section we will show, however,  that
it arises naturally as a degeneration of a corresponding Pieri formula for the Macdonald polynomials associated with root systems \cite{mac:orthogonal,mac:affine}.

If the bound $\langle \omega,\alpha^\vee\rangle= 2$ in Definition \ref{small:def} is never attained,
then the weight $\omega$ ($\neq 0$) is usually referred to as  {\em minuscule} and when the bound is attained only {\em once} it is called  {\em quasi-minuscule} \cite{bou:groupes}. 
Our Pieri formula then simplifies respectively to
\begin{subequations}
\begin{equation}
m_{\omega;R} P_\lambda =
\sum_{\substack{ \nu\in W \omega \\ \lambda +\nu \in P^+}}  V_\nu (\lambda ) P_{\lambda +\nu} \label{PFmin}  
\end{equation}
if $\omega$ is minuscule, and to
\begin{equation}
\bigl( m_{\omega;R}-m_{\omega;R} (\rho_{t,R^\vee}) \bigr) P_\lambda =
\sum_{\substack{ \nu\in W \omega \\ \lambda +\nu \in P^+}}  V_\nu (\lambda ) \bigl( P_{\lambda +\nu} -P_\lambda \bigr) \label{PFqmin}
\end{equation}
if $\omega$ is quasi-minuscule, where
\begin{equation}\label{Vmin}
V_\nu (\lambda ) :=
e^{-\nu}(\rho_{t,R^\vee})
\prod_{\substack{\alpha\in R^+_\lambda \\    \langle \alpha^\vee ,\nu\rangle > 0}} \frac{1-t_\alpha e^{\alpha}(\rho_{t,R^\vee}) }{1-e^{\alpha}(\rho_{t,R^\vee})  }  .
\end{equation}
\end{subequations}

By varying $\omega$ over the small {\em fundamental} weights,
$E_{\omega;R}$ \eqref{Eomega} produces various algebraically independent elements in
$\mathbb{C}[P]^W$, since the expansion on the monomial basis is of the triangular form
\begin{equation}\label{E-mon}
E_{\omega;R}=m_{\omega;R} + \sum_{\mu <\omega} \epsilon_{\omega\mu}(t) m_{\mu;R} ,  \quad
 \epsilon_{\omega\mu}(t) \in \mathbb{R}. 
\end{equation}
For the classical root systems all fundamental weights are small and we thus obtain a complete system of generators for $\mathbb{C}[P]^W$ this way. For the exceptional root systems some fundamental weights fail to be small so only part of such a generating system is recovered. Specifically,   the small fundamental weights read for the exceptional root systems $R=E_6$: $\omega_1,\omega_2,\omega_3,\omega_5,\omega_6$,
$R=E_7$:  $\omega_1,\omega_2,\omega_6,\omega_7$, $R=E_8$:  $\omega_1, \omega_8$, $R=F_4$: $\omega_1,\omega_4$,  and $R=G_2$: $\omega_1$, where we have numbered the fundamental weight in accordance with the tables in Bourbaki \cite{bou:groupes}. (Apart from
the small fundamental weights the class of small weights consists of the zero weight and all sums of two not necessarily distinct minuscule weights.)
For $R=A_n$ the fundamental weights are not only small but in fact minuscule and the corresponding expressions of the type in
Eqs. \eqref{PFmin}, \eqref{Vmin} reproduce the well-known Pieri formulas  for the Hall-Littlewood polynomials due to Morris \cite{mor:note}.

\begin{remark}\label{zeros:rem}
In the formula of Theorem \ref{PF:thm}  the coefficient of $P_{\lambda +\nu_1}$ vanishes when
$\lambda +\nu_1\not\in P^+$ due to a zero in $V_{\nu_1;-}(\lambda )$. Indeed, in this situation there exists a simple root $\beta\in R^+$ such that $\langle \lambda + \nu_1,\beta^\vee \rangle < 0$ (whence
$\langle \lambda ,\beta^\vee \rangle \leq 1$ by Lemma \ref{sw-chain:lem}).
One thus picks up a zero from the factor
$1-t_{\beta}e^{-\beta}(\rho_{t,R^\vee})$ in the numerator of  $V_{\nu_1;-}(\lambda )$.
\end{remark}

\begin{remark}
For $R=A_n$ the polynomials $E_{\omega ;R}$ with $\omega$ fundamental amount to the elementary symmetric functions. For $R=B_n$ the polynomials in question turn out to be particular instances of Okounkov's $BC$-type interpolation polynomials \cite[Appendix C] {kom-nou-shi:kernel}.
An alternative characterization elucidating the structure of the polynomials $E_{\omega ;R}$ for arbitrary $R$ in terms of Macdonald polynomials is given in Remark \ref{MPsw}  below.
\end{remark}

\section{Pieri formula for Macdonald polynomials}\label{sec3}
Let $S$ be either $R$ or $R^\vee$, i.e. $(R,S)$ constitutes an admissible pair of root systems in the sense of Ref. \cite{mac:orthogonal}, and  let $q_\alpha :=q^{u_\alpha}$ with $u_\alpha :=1$ if $S=R$ and $u_\alpha :=\langle\alpha ,\alpha\rangle /2$ if $S=R^\vee$.
For a formal series $f=\sum_{\lambda\in P} f_\lambda e^\lambda$, $f_\lambda\in \mathbb{C}$, we define $\int f := f_0$ and $\bar{f}:=\sum_{\lambda\in P}\bar{f}_\lambda e^{-\lambda}$ (with $\bar{f}_\lambda$ meaning the complex conjugate of $f_\lambda$).  
The Macdonald inner product on $\mathbb{C}[P]$ is then given by \cite{mac:orthogonal,mac:affine}
\begin{equation}\label{ipa}
\langle f,g\rangle_\Delta := |W|^{-1} \int  f \bar{g} \Delta \qquad (f,g\in\mathbb{C}[P]) ,
\end{equation}
with $|W|$ denoting the order of $W$ and
\begin{equation}\label{ipb}
\Delta := \prod_{\alpha \in R} \frac{( e^\alpha;q_\alpha)_\infty}{( t_\alpha e^\alpha;q_\alpha)_\infty} ,
\end{equation}
where we have employed the standard notation for the $q$-shifted factorial $(a;q)_m :=\prod_{k=0}^{m-1} (1-aq^k)$ with $m$ nonnegative integral or $\infty$.

\begin{definition}[Macdonald Polynomials \cite{mac:orthogonal,mac:affine}]
For $\lambda\in P^+$, the Macdonald polynomial is defined as the unique element in $\mathbb{C}[P]^W$ of the form
\begin{subequations}
\begin{equation}\label{mp-d1}
P_\lambda^{R,S} = \sum_{\mu \leq  \lambda } c_{\lambda \mu} (q,t)\, m_{\mu;R}
\end{equation}
with $c_{\lambda \mu} (q,t)\in\mathbb{C}$ such that
\begin{equation}
c_{\lambda \lambda}(q,t) = q^{\langle \lambda,\rho_{t,R^\vee} \rangle}
\prod_{\alpha\in R^+}\frac{(q_\alpha^{\langle \alpha^\vee , \rho_{t,S}\rangle };   q_\alpha)_{\langle \alpha^\vee ,\lambda\rangle}}{(t_\alpha q_\alpha^{\langle \alpha^\vee ,\rho_{t,S}\rangle};q_\alpha)_{\langle \alpha^\vee,\lambda\rangle}} 
\end{equation}
and
  \begin{equation}\label{mp-d2}
  \langle P_\lambda^{R,S} , m_{\mu ;R} \rangle_\Delta =0 \quad \text{for all} \ \mu < \lambda  
  \end{equation}
 \end{subequations}
 (where the normalization is chosen such that $P^{R,S}_\lambda (\rho_{t;R^\vee})=1$).
\end{definition}

The Pieri formula of Theorem \ref{PF:thm} arises as a degeneration of the following Pieri formula for the Macdonald polynomials.
\begin{theorem}[Pieri Formula for Macdonald Polynomials]\label{PFq:thm}
Let $\lambda, \omega \in P^+$ with $\omega $ small. Then
\begin{subequations}
\begin{equation}\label{PFq1}
E_{\omega;R} P_{\lambda}^{R,S} = 
\sum_{\substack{\mu_1<\cdots <\mu_\ell=\omega\\ \ell =1,\ldots ,n_\omega}}
(-1)^{\ell -1} \!\!\!\!\!\!\!\!
\sum_{\substack{\nu_k\in W_{\nu_{k-1}}  (  w^{-1}_{\nu_{k-1}} \mu_k  )   \\ k=1,\ldots ,\ell}} \!\!\!\! \!\!\!\!
 P_{\lambda +\nu_1}^{R,S}
\prod_{1\leq k\leq \ell}
V^{R,S}_{\nu_k;\nu_{k-1}}(\lambda +\rho_{t;S})  , 
\end{equation}
with $E_{\omega;R}$ taken from Eq. \eqref{Eomega} and
\begin{equation}\label{PFq2}
V_{\nu;\eta}^{R,S}(x) :=
\prod_{\substack{\alpha\in R_{\eta} \\     \langle \alpha^\vee ,\nu\rangle>0 }}
t_\alpha^{-\langle \alpha^\vee ,\nu \rangle/2 }
\frac{(t_\alpha q_\alpha^{\langle \alpha^\vee,x \rangle};q_\alpha )_{\langle \alpha^\vee ,\nu\rangle }}
{( q_\alpha^{\langle \alpha^\vee,x\rangle};q_\alpha )_{\langle \alpha^\vee ,\nu\rangle }} 
\end{equation}
\end{subequations}
(with the conventions $R_{\nu_0}=R$, $W_{\nu_0}=W$, $w_{\nu_0}=\text{Id}$ and 
$V^{R,S}_{\nu_1;\nu_{0}}=V^{R,S}_{\nu_1; 0}=V^{R,S}_{\nu_1;-}$).
\end{theorem}
In this Pieri formula  the coefficient of
$ P_{\lambda +\nu_1}^{R,S}$ again vanishes when
$\lambda +\nu_1\not\in P^+$, due to a zero in
$V^{R,S}_{\nu_1;-}(\lambda +\rho_{t;S} )$  stemming from the factor(s) in the numerator of the form
$(t_{\beta} q_{\beta}^{\langle -\beta^\vee,\lambda +\rho_{t;S} \rangle};q_{\beta} )_{\langle -\beta^\vee ,\nu_1\rangle } $ with $\beta\in R^+ $  simple such that $\langle \lambda + \nu_1,\beta^\vee \rangle < 0$
 (cf. Remark \ref{zeros:rem}).

It is well-known that for $q\to 0$ the Macdonald polynomial $P_\lambda^{R,S}$ tends to the Macdonald spherical function $P_\lambda$ (so the dependence on the choice of $S$ disappears in the limit) \cite{mac:orthogonal}. Performing this limit transition
at the level of the Pieri formulas readily entails Theorem \ref{PF:thm} as a consequence of Theorem \ref{PFq:thm}.

\begin{remark}
For $\omega$ (quasi-)minuscule (and $S=R$) the Pieri formula \eqref{PFq1}, \eqref{PFq2} can be found for instance in Ref. \cite[Appendix A]{die:scattering} and---with a more detailed derivation---in Ref. \cite{las:inverse}.
When $R=A_n$ and $\omega$ is fundamental, the formula in question reproduces the well-known Pieri formulas for the Macdonald symmetric functions \cite{mac:symmetric}. For other classical root systems related Pieri formulas for the corresponding Macdonald polynomials
can be found in Refs. \cite{die:self-dual,las:inverse} (cf. also Ref. \cite{die-ito:difference} for the case of the exceptional root system $G_2$). Pieri-type formulas for the nonsymmetric Macdonald polynomials were recently studied in Ref. \cite{bar:pieri-type}.
\end{remark}

\begin{remark}\label{MPsw}
The Macdonald polynomial associated with a small weight $\omega$ is easily computed explicitly by means of the Pieri formula. Indeed, setting $\lambda=0$ in Theorem \ref{PFq:thm} yields
\begin{subequations}
\begin{equation}\label{inverse}
E_{\omega;R}  = 
\sum_{\substack{\mu_1<\cdots <\mu_\ell=\omega\\ \ell =1,\ldots ,n_\omega}}
(-1)^{\ell -1}  p_{\mu_1}^{R,S}
\prod_{1\leq k\leq \ell-1}
V^{R,S}_{\mu_{k+1};\mu_{k}}(\rho_{t;S})  ,
\end{equation}
where
$p_{\mu}^{R,S}:= V^{R,S}_{\mu ;-}(\rho_{t;S}) P_{\mu}^{R,S}=c_{\mu\mu}^{-1}(q,t)P_{\mu}^{R,S}$.
This relation is the inverse of the following unitriangular 
expansion of the (monic) Macdonald polynomial $p_{\omega}^{R,S}$ in terms of the generators $E_{\mu;R}  $, $\mu\in P^+_\omega$:
\begin{equation}\label{explicit}
 p_{\omega}^{R,S}
= E_{\omega;R}+\sum_{\mu< \omega}
E_{\mu;R}  V^{R,S}_{\omega;\mu}(\rho_{t;S})  .
\end{equation}
\end{subequations}
The obtained expression in Eq. \eqref{explicit} is a reduced root system counterpart of the
expansion in 
Ref \cite[Theorem 5.1]{kom-nou-shi:kernel} for the case of the nonreduced root systems. Reversely, the formula in question characterizes our polynomials $E_{\omega;R}$ as the unique polynomials of the form in Eq.  \eqref{E-mon} on which
the monic Macdonald polynomials $p_{\omega}^{R,S}$ expand triangularly
with simple explicit factorized coefficients of the form $ V^{R,S}_{\omega;\mu}(\rho_{t;S})$.
It is instructive to notice that for $t=1$  the coefficients $V^{R,S}_{\mu_{k+1};\mu_{k}}(\rho_{t;S}) $ are equal to $1$ and
$p_{\mu}^{R,S}=m_{\mu;R}$. In this case Eq. \eqref{inverse} amounts to Hall's identity for the
M\"obius inversion of Eq. \eqref{explicit}. The case of general $t$ boils down to Gauss' formula for the inverse of a triangular matrix.
\end{remark}

\begin{remark}
By combining our Pieri formula with the representation of the
monic Macdonald polynomials associated with the small weights in
Eq. \eqref{explicit},
one arrives at a special Littlewood-Richardson-type rule (or linearization formula) for the products of two Macdonald polynomials with one of the weights being small:
\begin{eqnarray}\label{PFqalt}
\lefteqn{P_{\omega}^{R,S}P_{\lambda}^{R,S}=} && \\
&& \!\!\!\!
c_{\omega\omega} (q,t)\sum_{\substack{\mu_1<\cdots <\mu_\ell\leq \omega\\ \ell =1,\ldots ,n_\omega}}
(-1)^{\ell -1}V^{R,S}_{\omega; \mu_\ell}(\rho_{t;S})  \!\!\!\!\!\!\!\! 
\sum_{\substack{\nu_k\in W_{\nu_{k-1}}  (  w^{-1}_{\nu_{k-1}} \mu_k  )   \\ k=1,\ldots ,\ell}} \!\!\!\! \!\!\!\!
 P_{\lambda +\nu_1}^{R,S}
\prod_{1\leq k\leq \ell}
V^{R,S}_{\nu_k;\nu_{k-1}}(\lambda +\rho_{t;S})  . \nonumber
\end{eqnarray}
The latter formula reveals that the terms in the expansion
$$P_{\omega}^{R,S}P_{\lambda}^{R,S}=
\sum_{\mu\leq\omega} \sum_{\substack{ \nu\in W\mu \\ \lambda+\nu\in P^+}}  
C_{\omega,\lambda }^{R,S}(\nu) P^{R,S}_{\lambda +\nu}$$
have completely factorized (structure) coefficients $C_{\omega,\lambda }^{R,S}(\nu) =  c_{\omega\omega}(q,t)V^{R,S}_{\nu; -}(\lambda +\rho_{t;S})$ for $\nu$ in the highest-weight orbit $W\omega$, whereas the coefficients
$C_{\omega,\lambda }^{R,S}(\nu) $, $\nu\in W\mu$ with $\mu < \omega$ 
are represented in contrast by (large) alternating sums of products. 
\end{remark}

\section{Difference equation for Macdonald polynomials}\label{sec4}
Given the admissible pair $(R,S)$, let us associate with $\alpha \in R$ the rescaled root
$\alpha_*:=u_\alpha^{-1}\alpha\in S$.  Furthermore, we will employ the notation
$P^+(R)$ and $P^+(R^\vee)$ to distinguish the cone of dominant weights associated with $R$ and $R^\vee$ (endowed with the corresponding dominance ordering stemming from $R^+$ and
$(R^\vee)^+$, respectively).

\begin{definition}[Generalized Macdonald Operator]
For $\omega\in P^+(S^\vee)$ small, the generalized Macdonald operator is given by
\begin{subequations}
\begin{equation}\label{macop1}
D_{\omega}^{R,S} := 
\sum_{\substack{\mu_1<\cdots <\mu_\ell=\omega\\ \ell =1,\ldots ,n_\omega}}
(-1)^{\ell -1} \!\!\!\!\!\!\!\!
\sum_{\substack{\nu_k\in W_{\nu_{k-1}}  (  w^{-1}_{\nu_{k-1}} \mu_k  )   \\ k=1,\ldots ,\ell}}
\Bigl( 
\prod_{1\leq k\leq \ell}
\hat{V}^{R,S}_{\nu_k;\nu_{k-1}} 
\Bigr)   \,
T_{\nu_1}  ,
\end{equation}
where 
\begin{equation}\label{macop2}
\hat{V}_{\nu;\eta}^{R,S} :=
\prod_{\substack{\alpha\in R_{\eta} \\     \langle \alpha_* ,\nu \rangle>0 }}
t_\alpha^{-\langle \alpha_* ,\nu \rangle/2 }
\frac{(t_\alpha e^\alpha;q_\alpha )_{\langle \alpha_* ,\nu \rangle }}
{( e^\alpha ;q_\alpha )_{\langle \alpha_* ,\nu \rangle }} 
\end{equation}
and $T_x:\mathbb{C}[P]\to\mathbb{C}[P]$, $x\in E$, denotes the translation operator
determined by $T_x e^\lambda := q^{\langle\lambda ,x\rangle}e^\lambda$, $\lambda\in P$.
\end{subequations}
\end{definition}
For nonreduced root systems analogous difference
operators were introduced in Ref. \cite{die:commuting}.
When $\omega$ is minuscule or quasi-minuscule, the operator
$D_{\omega}^{R,S}$ reduces to the Macdonald difference operators \cite{mac:orthogonal}
\begin{subequations}
\begin{equation}\label{macop:m}
D_\omega^{R,S}=q^{-\langle \omega,\rho_{t;S}\rangle} \sum_{\nu\in W \omega} 
\Bigl( \prod_{\substack{\alpha\in R \\     \langle \alpha_* ,\nu \rangle>0 }}
\frac{1-t_\alpha e^\alpha }
{1-e^\alpha } \Bigr) \, T_\nu  
\end{equation}
and
\begin{equation}\label{macop:qm}
D_\omega^{R,S}=q^{-\langle \omega,\rho_{t;S}\rangle} \sum_{\nu\in W \omega} 
\Bigl(
\prod_{\substack{\alpha\in R \\     \langle \alpha_* ,\nu \rangle>0 }}
\frac{(t_\alpha e^\alpha;q_\alpha )_{\langle \alpha_* ,\nu \rangle }}
{( e^\alpha ;q_\alpha )_{\langle \alpha_* ,\nu \rangle }} \Bigr)
\, \bigl( T_\nu -1\bigr) ,
\end{equation}
\end{subequations}
respectively.
The following theorem generalizes the corresponding diagonal action of $D_\omega^{R,S}$ \eqref{macop:m}, \eqref{macop:qm}
on the Macdonald basis  \cite{mac:orthogonal} to the case that $\omega$ is an arbitrary small weight in $P^+(S^\vee)$.

\begin{theorem}[Difference Equation]\label{Deq:thm}
Let $\omega\in P^+(S^\vee)$ be small. Then for any $\lambda\in P^+(R)$
\begin{equation}
D_\omega^{R,S}  P_\lambda^{R,S}=
E_{\omega;S^\vee} (\lambda +\rho_{t;S})  P_\lambda^{R,S} .
\end{equation}
\end{theorem}
The Pieri formula in Theorem \ref{PFq:thm}  follows from the difference equation in Theorem
\ref{Deq:thm} by a standard argument (see e.g. Refs. \cite{mac:symmetric,mac:affine,die:self-dual,die-ito:difference,las:inverse}) involving
the {\em Duality Symmetry}  \cite[Sec. 5.3]{mac:affine}
 \begin{equation}\label{duality:eq}
P_{\lambda}^{R,S}(\mu+\rho_{t;R^\vee})=
   P_{\mu}^{S^\vee,R^\vee}(\lambda+\rho_{t;S})
   \end{equation}
for all $\lambda\in P^+(R),\mu\in P^+(S^\vee )$. Indeed, evaluating the Pieri formula \eqref{PFq1}, \eqref{PFq2} at $\mu +\rho_{t;R^\vee}$, $\mu\in P^+(S^\vee)$ and
applying the duality symmetry to both sides entails the difference equation
$D_\omega^{S^\vee,R^\vee}  P_\mu^{S^\vee,R^\vee}=
E_{\omega;R} (\mu +\rho_{t;R^\vee})  P_\mu^{S^\vee,R^\vee}$ evaluated at $x=\lambda+\rho_{t;S}$ (since $V^{R,S}_{\nu ;\eta} (x)=\hat{V}^{S^\vee ,R^\vee} _{\nu ; \eta} (x)$). Invoking of Theorem \ref{Deq:thm} thus verifies the Pieri formula, first for $x=\mu +\rho_{t;R^\vee}$, $\mu\in P^+(S^\vee)$ 
and then as an identity in $\mathbb{C}[P]^W$
by the polynomiality of both sides.

\begin{remark}\label{W-invariance:rem}
The generalized Macdonald operator is $W$-invariant:
$w D_\omega^{R,S}w^{-1} = D_\omega^{R,S}$ for all $w\in W$. This invariance hinges on the relations
$w\hat{V}_{\nu;\eta}^{R,S} =\hat{V}_{w\nu;w\eta}^{R,S} $ and
$wT_\nu w^{-1} =T_{w\nu}$ upon invoking the stability properties of small weights
in Lemma \ref{sw-stab:lem}.
\end{remark}

\section{Proof of the difference equation}\label{sec5}

The proof of the difference equation in Theorem \ref{Deq:thm} is based on the following three properties of the generalized Macdonald operator 
$D^{R,S}_\omega $ \eqref{macop1}, \eqref{macop2}  with $\omega\in P^+(S^\vee)$ small.

\begin{proposition}\label{residue:prp}
The operator $D^{R,S}_\omega $ maps $\mathbb{C}[P]^W$ into itself.
\end{proposition}

\begin{proposition}\label{triangularity:prp}
The action of $D^{R,S}_\omega $ in  $\mathbb{C}[P]^W$
is triangular with respect to the monomial basis:
$$D_\omega^{R,S} m_{\lambda ;R} = E_{\omega ;S^\vee} (\lambda+\rho_{t,S}) m_{\lambda ;R}+\sum_{ \mu < \lambda} b_{\lambda\mu}(q,t)m_{\mu ;R},\qquad
\forall \lambda\in P^+(R),$$ 
with $b_{\lambda\mu}(q,t)\in\mathbb{C}$.
\end{proposition}

\begin{proposition}\label{symmetry:prp}
The operator $D^{R,S}_\omega $ is symmetric with respect to Macdonald's inner product $\langle\cdot,\cdot \rangle_\Delta$ \eqref{ipa}, \eqref{ipb}:
$$\langle D_\omega^{R,S} f ,g\rangle_\Delta=\langle  f , D_{\omega}^{R,S} g\rangle_\Delta ,
\qquad\forall  f,g\in \mathbb{C}[P]^W.$$
\end{proposition}

It is evident from Propositions \ref{residue:prp}--\ref{symmetry:prp} that for generic $\epsilon\in\mathbb{R}$ the expression
$(E_{\omega ;S^\vee} (\lambda+\rho_{t,S}) +\epsilon)^{-1} (D_\omega^{R,S} +\epsilon )P_\lambda^{R,S}$ satisfies
the defining properties in Eqs. \eqref{mp-d1}--\eqref{mp-d2} of the Macdonald polynomial $P_\lambda^{R,S}$, which proves the difference equation.

\subsection{Proof of Proposition \ref{residue:prp}}
Let $f\in\mathbb{C}[P]^W$.
To infer that $D^{R,S}_\omega f\in \mathbb{C}[P]^W$ it is enough to show
that  $D^{R,S}_\omega f\in \mathbb{C}[P]$ in view of Remark \ref{W-invariance:rem} above.
To this end we will check that the evaluated expression
$(D^{R,S}_\omega f) (x)$ is regular as a function of $x$ in (the complexification of) $E$.
The expression in question is obtained from 
$$
\sum_{\substack{\mu_1<\cdots <\mu_\ell=\omega\\ \ell =1,\ldots ,n_\omega}}
(-1)^{\ell -1} 
f (x+\mu_1)
\prod_{1\leq k\leq \ell}
\hat{V}^{R,S}_{\mu_k;\mu_{k-1}} (x) 
$$
by symmetrization in $x$ with respect to the action of the Weyl group.
In the term associated with a fixed ascending chain $\mu_1<\cdots <\mu_\ell=\omega$ 
the factor $\hat{V}^{R,S}_{\mu_k;\mu_{k-1}} $, $1\leq k \leq \ell$,
gives rise to poles caused by its
factors in the denominator of the form
$(e^\alpha ;q_\alpha)_{\langle \alpha_*,\mu_k\rangle}$,
$\alpha\in R^+_{\mu_{k-1}}\setminus R^+_{\mu_k}$.
The poles stemming from the factors of the form $(1-e^\alpha )$ manifestly cancel
in $(D^{R,S}_\omega f) (x)$ due to the Weyl-group invariance.  Poles caused by factors of the form $(1-q_\alpha e^\alpha )$, however, require a more careful analysis.
The latter poles occur for  $\alpha\in R^+_{\mu_{k-1}}$
when $\langle \alpha_*,\mu_k\rangle=2$.
In this situation we may exploit the $W$-invariance of $(D^{R,S}_\omega f) (x)$
and assume that $\mu_k>\mu_k-\alpha_*^\vee \geq\mu_{k-1}$ in view of
Lemma \ref{sw-cover:lem}  below.
We will now show that the pole of interest cancels pairwise against a corresponding pole
in the term associated with the chain obtained from $\mu_1<\cdots <\mu_\ell=\omega$
by inserting or extracting  $\mu_k-\alpha_*^\vee$ into or from the chain
depending on whether $\mu_k-\alpha_*^\vee>\mu_{k-1}$ or
$\mu_k-\alpha_*^\vee=\mu_{k-1}$, respectively (with the convention that we perform an insertion if $k=1$).
Indeed, by Property (i) of Lemma \ref{sw-stab:lem} it is clear that
$\alpha\in R^+_{\mu_{k-2}}$, which ensures that the presence of the pole is stable with respect to
such insertion/extraction operations.
Furthermore, since  these operations on the chain determine an involution
we may assume without restriction  
that we are dealing with the (insertion) case $\mu_k-\alpha_*^\vee>\mu_{k-1}$. 
After dividing out the common factor
$ \prod_{\substack{1\leq j\leq \ell \\ j\neq k}}
\hat{V}^{R,S}_{\mu_j;\mu_{j-1}} $,  we see that the residue of the generically simple poles
at the hyperplane $\langle \alpha_*,x\rangle =-1$ in the terms associated with
the chains $\mu_1<\cdots <\mu_{k-1}<\mu_k<\cdots  <\mu_\ell=\omega$ and
$\mu_1<\cdots <\mu_{k-1}<\mu_k-\alpha_*^\vee<\mu_k<\cdots  <\mu_\ell=\omega$
cancel pairwise, since
$$ \text{Res}_{\langle \alpha_*,x\rangle =-1}  \left[ \hat{V}^{R,S}_{\mu_k;\mu_{k-1}} (x) \right]
= \text{Res}_{\langle \alpha_*,x\rangle =-1} \left[ \hat{V}^{R,S}_{\mu_k;\mu_{k}-\alpha_*^\vee} (x) \hat{V}^{R,S}_{\mu_k-\alpha_*^\vee;\mu_{k-1}} (x) \right]  
$$
and if $k=1$ moreover
$$f(x+\mu_1)\bigr|_{\langle \alpha_*,x\rangle =-1} =
f(x+\mu_1-\alpha_*^\vee )\bigr|_{\langle \alpha_*,x\rangle =-1}  .$$
The equality of the evaluations is immediate from the $W$-invariance of $f$ because $x+\mu_1-\alpha_*^\vee=r_\alpha (x+\mu_1)$ for $x$ in the hyperplane $\langle \alpha_*,x\rangle =-1$
(where $r_\alpha: E\to E$ refers to the orthogonal reflection across the hyperplane perpendicular to $\alpha$).
The equality of the residues follows from a case-by-case check, varying $\mu_k$ over the list of small weights and picking $\alpha$ as in Part (ii) of Lemma \ref{sw-cover:lem}  (taking---without sacrificing generality---$R_{\mu_{k-1}}=R$).

The upshot is that in the $W$-invariant expression $(D^{R,S}_\omega f) (x)$ all poles cancel out to produce a regular expression in $x$.

\subsection{Proof of Proposition \ref{triangularity:prp}}
By Proposition \ref{residue:prp},  $D^{R,S}_\omega m_{\lambda ;R}$
expands as a finite linear combination of monomials $m_{\mu ;R}$, $\mu\in P^+(R)$.
For $x$ growing to infinity in such a way that
$\langle x, \alpha \rangle\to -\infty$ for all $\alpha\in R^+$, clearly
$\hat{V}_{\nu;\eta}^{R,S} (x)=q^{\langle \nu, \rho_{t;S_{\eta}} \rangle}(1+o(1))$ and $m_{\lambda;R}(x)= q^{\langle \lambda , x \rangle}(1 + o(1))$, which entails that
$$(D^{R,S}_\omega m_{\lambda ;R}) (x)= E_{\omega;S^\vee} (\lambda +\rho_{t;S})
 q^{\langle \lambda , x \rangle}\bigl( 1 + o(1)\bigr) $$
(since
$\sum_{\nu\in W_\eta (w_\eta^{-1} \mu )} q^{\langle \nu, \rho_{t;S_{\eta}} \rangle}=
m_{\mu ; S^\vee_{w_\eta \eta}}(\rho_{t;S_{w_\eta \eta}})$).
From this asymptotics deep in the anti-dominant Weyl chamber it is immediate that the monomial expansion of $D^{R,S}_\omega m_{\lambda ;R}$ is of the form stated in the proposition.

\subsection{Proof of Proposition \ref{symmetry:prp}}
Let us momentarily assume that $t_\alpha$ is an integral power $\geq 2$ of $q_\alpha$. This parameter restriction guarantees that all terms of $\Delta D^{R,S}_\omega f$ lie in $\mathbb{C}[P]$.
The proposition is now immediate (upon symmetrization with respect to the action of the Weyl group) from the following adjointness relations associated to any
chain of small weights in $P^+(S^\vee)$ of the form $\mu_1<\mu_2<\cdots <\mu_\ell$:
$$
\langle  \hat{V}^{R,S}_{\nu_{1};-} \prod_{1\leq k\leq \ell-1}
\hat{V}^{R,S}_{\nu_{k+1};\nu_{k}}
T_{\nu_1} f,g\rangle_\Delta =
 \langle f,  \hat{V}^{R,S}_{\nu_{1};-} \prod_{1\leq k\leq \ell-1}
\hat{V}^{R,S}_{w_{0}(R_{\nu_k}) \nu_{k+1};\nu_{k}} 
 T_{\nu_1}g\rangle_\Delta ,
$$
where $\nu_k\in W_{\nu_{k-1}}  (  w^{-1}_{\nu_{k-1}} \mu_k  ) $, $k=1,\ldots ,\ell$ (so, in particular,
the weights $\nu_1,\ldots,\nu_\ell$ all lie in the same Weyl chamber) and
$w_0(R_{\eta})$ denotes the longest element in the Weyl group $W_\eta$ of $R_\eta$.
These adjointness relations follow in turn from the elementary properties
$\int T_x (f)\overline{g}=\int f \overline{T_x (g)}$ for $x\in E$, 
$\overline{\hat{V}^{R,S}_{\nu;\eta} }= \hat{V}^{R,S}_{-\nu;\eta}=\hat{V}^{R,S}_{w_{0}(R_{\eta})\nu;\eta}$ and
$T_\nu \hat{V}^{R,S}_{-\nu;-} \Delta=\hat{V}^{R,S}_{\nu;-} \Delta$ for $\nu$ in the Weyl orbit of a small weight,
together with the observation that the factors $\hat{V}^{R,S}_{\nu_{k+1};\nu_{k}} $, $k=1,\ldots ,\ell-1$, commute with $T_{\nu_1}$ (by Property (i) of Lemma \ref{sw-stab:lem}).

Finally, our temporary technical restriction on the values of multiplicity function $t$ is removed by exploiting that the Macdonald polynomials---and thus both sides of the eigenvalue equation in Theorem \ref{Deq:thm}---are rational expressions in $t_\alpha$
(which allows to remove this restriction from Theorem \ref{Deq:thm})
and that the Macdonald polynomials moreover form an orthogonal basis of $\mathbb{C}[P]^W$ with
inner product $\langle \cdot ,\cdot\rangle_\Delta$ (by which one subsequently concludes  that Proposition \ref{symmetry:prp} also holds for general $t$ since   $D^{R,S}_\omega $ is thus unitarily equivalent to a real multiplication operator)  \cite{mac:orthogonal,mac:affine}.

\vspace{3ex}

{\bf Acknowledgments.} We thank the referee for some helpful suggestions improving our presentation.

\appendix

\section{Key properties of small weights}\label{sw:app}
The following properties
are readily seen from a straightforward case-by-case analysis upon listing the small weights for all irreducible reduced root systems.

\begin{lemma}\label{sw-chain:lem}
For $\omega$ small, the highest-weight system
$P^+_\omega$ constitutes a linear chain of small weights.
\end{lemma}

\begin{lemma}\label{sw-stab:lem}
For any pair of small weights
$\mu_a<\mu_b$:
\begin{itemize}
\item[(i)] $\alpha\in R_{\mu_a}\setminus  R_{\mu_b}\Rightarrow\alpha\in  \bigcap_{\mu\leq\mu_a}R_{\mu},$  \vspace{1ex}
\item[(ii)] $ (\bigcap_{\mu\leq\mu_a}W_{\mu})(\mu_{b})=W_{\mu_{a}}(\mu_{b}).$ 
\end{itemize}
\end{lemma}

\begin{lemma}\label{sw-cover:lem}
If $\mu_a<\mu_b$ is a pair of small weights with $\beta\in R^+_{\mu_a}$ such that $\langle \mu_b , \beta^\vee \rangle =2$,
then
\begin{itemize}
\item[ (i)] $(W_{\mu_a } \cap W_{\mu_b } )(\beta)=
\{ \alpha\in R^+_{\mu_a} \mid \langle \mu_b ,\alpha^\vee \rangle =2, \langle\alpha ,\alpha\rangle = \langle\beta ,\beta\rangle \} ,$ \vspace{1ex}
\item[(ii)] moreover, there exists a unique root $\alpha\in (W_{\mu_a } \cap W_{\mu_b }) (\beta)$ such that $\mu_b -\alpha \geq \mu_a $.
\end{itemize}
\end{lemma}

\begin{remark}
(i) In  Lemma \ref{sw-chain:lem} is immediate
that any $\mu\in P^+_\omega$ is a small weight as the highest root $\alpha_0^\vee$
of $R^\vee$ is situated in the dominant Weyl chamber $\{ x\in E\mid \langle x, \alpha \rangle \geq 0,\ \forall \alpha\in R^+\}$, so $\langle \mu, \alpha^\vee\rangle\leq
\langle \mu, \alpha_0^\vee\rangle\leq \langle \omega, \alpha_0^\vee\rangle \leq 2$
for all $\alpha\in R^+$.  

(ii) In Lemma \ref{sw-stab:lem} Part (ii) follows from Part (i).
Indeed, since
$\bigcap_{\mu\leq\mu_a}W_{\mu}\subset W_{\mu_a}$ are the Weyl groups of the root systems
$\bigcap_{\mu\leq\mu_a}R_{\mu}\subset R_{\mu_a}$, respectively,
a count of the size of the Weyl orbits on both sizes of Part (ii) with the aid of  the general formulas
$|W(\mu)|= |W|/|W_\mu |$ and $|W|=\prod_{\alpha \in R^+} (1+\text{ht}(\alpha))/\text{ht}(\alpha )$ (where $\text{ht}(\alpha)$ refers to the height of $\alpha$ with respect to the basis of simple roots) 
reveals that Part (i) implies that (the size of) both orbits must be equal.

(iii)  When checking  Lemma \ref{sw-cover:lem} it is enough to consider the case $\mu_a=0$.
Furthermore, in Part (ii) the weight $\mu_b-\alpha$ 
is the predecessor of $\mu_b$ if (and only if) $\alpha$ is short (cf. \cite[Theorem 2.6]{ste:partial}).
\end{remark}

\bibliographystyle{amsplain}

\end{document}